\documentclass[11pt]{amsart}
\usepackage{mathrsfs}
\usepackage{mathpple}
\usepackage{amsmath}
\usepackage{amssymb}
\usepackage{mathtools}
\usepackage{bm}
\usepackage{array}
\usepackage{booktabs}
\usepackage{needspace}
\usepackage[normalem]{ulem}
\usepackage{hyperref,xcolor}
\definecolor{wine-stain}{rgb}{0.5,0,0}
\hypersetup{
	colorlinks,
	linkcolor=wine-stain,
	citecolor=wine-stain,
	urlcolor=wine-stain,
	linktoc=all
}

\newtheorem{theorem}{Theorem~}[section]
\newtheorem{lemma}[theorem]{Lemma~}
\newtheorem{proposition}[theorem]{Proposition~}
\newtheorem{corollary}[theorem]{Corollary~}

\theoremstyle{remark}
\newtheorem{remark}[theorem]{Remark~}

\theoremstyle{definition}

\newcommand{\abs}[1]{\left\vert#1\right\vert}
\newcommand{\BF}{\mathbb F}
\newcommand{\BQ}{\mathbb Q}
\newcommand{\BZ}{\mathbb Z}
\newcommand{\Aut}{\operatorname{Aut}}
\newcommand{\Co}{\operatorname{Co}}
\newcommand{\NS}{\operatorname{NS}}

\newcommand{\Syl}{\operatorname{Syl}}
\newcommand{\rk}{\operatorname{rank}}

\begin{document}

\title[Finite Symplectic Groups in Characteristic Seven]
{Finite Symplectic Automorphism Groups of K3 Surfaces in
Characteristic Seven}
\author{Bin Wang, Zhiwei Zheng}
\address{}
\email{}

\newcommand{\Addresses}{{%
		\bigskip
		\footnotesize

		B.~Wang, \textsc{Hetao Institute of Mathematics and Interdisciplinary Sciences, Shenzhen, China}\par\nopagebreak
		\textit{Email address}: \texttt{matph92@gmail.com}

		\medskip

		Z.~Zheng, \textsc{Tsinghua University, Beijing, China}\par\nopagebreak
		\textit{Email address}: \texttt{zhengzhiwei@mail.tsinghua.edu.cn}
}}

\thanks{}

\begin{abstract}
We classify all finite groups of order divisible by $7$ that admit a faithful symplectic action on a K3 surface over an algebraically closed field of characteristic seven, together with the possible heights and Artin invariants.  Combined with the known classification of tame actions, this gives the complete list of finite groups admitting such an action in characteristic seven.
\end{abstract}
\maketitle


\section{Introduction}

Mukai classified the finite groups that act faithfully and
symplectically on complex K3 surfaces \cite{Muk88}.  Over an algebraically closed field of positive characteristic,  the classification
is known in tame case \cite{DK09}, \cite{OS24}.  Dolgachev--Keum
\cite[Theorem 1.1]{DK09J} gave the classification of wild actions in
characteristic $11$. 

In this paper, we focus on K3 surfaces in characteristic $7$. We call a finite group $G$ \emph{wild} or \emph{tame} according as $7\mid\abs G$ or $7\nmid\abs G$. We will classify the wild case.  

Write $\Aut^s(X)$ for the subgroup of automorphisms acting trivially on $H^0(X,\Omega_X^2)$. For
$G\leq\Aut^s(X)$ finite, put $S_G(X)=(\NS(X)^G)^\perp\subset\NS(X)$.
When $X$ is supersingular, the saturation of $(G,S_G(X))$ is
$(\widetilde G,S_G(X))$, where $\widetilde G$ is the full group of
isometries of $S_G(X)$ acting trivially on its discriminant group.
Put $L_2(7)=\operatorname{PSL}_2(\BF_7)$, and let
$F_{21}=C_7\rtimes C_3$ denote the nonabelian group of order $21$.

Our main result is the following classification.

\begin{theorem}
\label{thm:all-wild-seven}
Let $X$ be a K3 surface over an algebraically closed field of characteristic $7$, and let $G\leq\Aut^s(X)$ be finite with
$7\mid\abs G$.  Then exactly one of the following occurs.
\begin{enumerate}
\item The surface $X$ has finite height $6$, and $G\cong C_7$ or $C_{14}$.
\item The surface $X$ is supersingular of Artin invariant $1$, and $G$ is isomorphic to a subgroup of $M_{21}.2_2$ or $2^4{:}A_7$ whose order is divisible by $7$ \cite[Theorem 1.2(4) and Corollary 1.3]{OS24}.
\item The surface $X$ is supersingular of Artin invariant $2$, and
$G$ is isomorphic to one of
\[
C_7,\ F_{21},\ L_2(7),\ C_{14},\ C_2\times F_{21},\
C_2\times L_2(7).
\]
\item The surface $X$ is supersingular of Artin invariant $3$, and $G\cong C_7$ or $F_{21}$.
\end{enumerate}
Conversely, every listed possibility occurs.  In the finite-height case both groups occur on the same explicit height-six surface. In particular, if $X$ is supersingular with Artin invariant
$\sigma\geq4$, then every finite-order automorphism of $X$ is tame.
\end{theorem}

We also give the refined supersingular classification in terms of
maximal abstract groups.

\begin{proposition}
\label{prop:all-supersingular-seven}
Let $X$ be a supersingular K3 surface over an algebraically closed
field of characteristic $7$, with Artin invariant $\sigma$.
A finite group $G$ acts faithfully and symplectically on $X$ if and
only if it is isomorphic to a subgroup of one of the groups in the
corresponding row:
\[
\begin{array}{c|l}
\sigma&\text{maximal groups}\\\hline
1&M_{21}.2_2,\ 2^4{:}A_7\\
2&4^2.A_4,\ \Gamma_{25}a_1,\ C_3^2{:}C_4,\
\mathrm{SD}_{16},\ D_{10},\ C_2\times L_2(7)\\
3&C_2^4,\ D_8,\ S_3,\ F_{21}\\
4&C_2^3,\ C_4,\ C_3\\
5,6&C_2\\
7,\ldots,10&1.
\end{array}
\]
Every listed group, and every subgroup of it, acts on every surface
in the corresponding Artin stratum.  The row $\sigma=1$ is due to
Ohashi--Sch\"utt \cite[Theorem 1.2(4)]{OS24}.
\end{proposition}

Here $D_{2n}$ denotes the dihedral group of order $2n$, and
$\mathrm{SD}_{16}$ is the semidihedral group of order sixteen.
The groups $4^2.A_4$, $\Gamma_{25}a_1$, and $C_3^2{:}C_4$ are
H\"ohn--Mason nos.~37, 43, and 46, respectively, with SmallGroups
identifiers $(192,1023)$, $(64,138)$, and $(36,9)$
\cite[Table 1]{HM16}, \cite[Table 10.2]{Has12}.
We prove Proposition \ref{prop:all-supersingular-seven} at the end of
Section \ref{sec:fine-classification}.


The ingredients of the proof separate naturally according to height. Keum's description of an order-$7$ automorphism and its rational quotient \cite[Propositions 2.2 and 8.1]{Keu16}, combined with the slope decomposition, forces a finite-height surface to have formal Brauer height 6.  The finite-height lattice bound of Wang--Zheng \cite[Proposition 1.3]{WZroot} and the action on the crystalline transcendental part then reduce the group to $C_7$ or $C_{14}$.  

For supersingular surfaces, Leech lattice methods and character
averaging show that Sylow $7$-subgroups have order at most $7$, while
discriminant gluing gives $\sigma\leq3$ for wild actions.  For
$\sigma\geq2$, the wild case reduces to three H\"ohn--Mason saturations.

We also give an explicit elliptic family with a symplectic
$F_{21}$-action and generic Artin invariant two in Appendix
\ref{app:explicit-seven}.

\medskip

\noindent\textbf{AI Disclosure:} 
This research is a natural continuation of our previous work \cite{WZroot}. The initial framework and methods were formulated by the authors. In subsequent stages, LLM-based AI tools were used in exploring proof strategies, refining proofs and editing the text. All contents were verified and finalized by the authors, who take full responsibility.

\section{Supersingular K3 Lattices and Leech Embeddings}
\label{sec:seven-framework}

For an integral lattice $L$, write $L^\vee$ for its dual lattice and $A_L=L^\vee/L$ for its discriminant group.  If $A$ is a finite abelian group, let $\ell(A)$ denote its minimal number of generators.  For a prime $p$, let $A_p$ be the $p$-primary part of $A$ and put $\ell_p(A)=\ell(A_p)$; thus $\ell(A)=\max_p\ell_p(A)$.  For a finite quadratic form $q$ on $A$, set $q_p=q|_{A_p}$, and write $A[p]=\{a\in A:pa=0\}$.

Let $X$ be a supersingular K3 surface over an algebraically closed field of characteristic $7$.  Its N\'eron--Severi lattice has signature $(1,21)$ and discriminant group $A_{\NS(X)}\cong(\BZ/7\BZ)^{2\sigma}$, where $1\leq\sigma\leq10$ is the Artin invariant.  We write $N_{7,\sigma}$ for this lattice.  For a finite group $G$ acting on a lattice $L$, put
$S_G(L)=(L^G)^\perp_L$, i.e. the orthogonal complement of the invariant sublattice.

We recall the lattice-theoretic results of crystalline Torelli which will be used below.  The marked theorem is due to Ogus \cite[Theorem II$''$]{Ogu83}.  The discriminant criterion for symplectic automorphisms is \cite[Proposition 3.1]{Jan16}, and the resulting lattice-theoretic formulation is in \cite[\S 3.3]{WZroot}.  Namely, if $G\leq\Aut^s(X)$ is finite, then the action on $\NS(X)$ is faithful and trivial on $A_{\NS(X)}$, while $S_G(X)$ is negative-definite and rootless and $G$ acts trivially on $A_{S_G(X)}$.  The faithfulness also follows from \cite[Corollary 2.5]{Ogu79}.

Call $(G,S)$ a \emph{Leech pair} if $G$ acts faithfully on the rootless negative-definite even lattice $S$, one has $S^G=0$, and $G$ acts trivially on $A_S$. Its \emph{saturation} is $\widetilde G=\ker\bigl(\Aut(S)\longrightarrow\Aut(A_S)\bigr)$. The pair is saturated when $G=\widetilde G$.  H\"ohn--Mason classify
the saturated Leech pairs arising from fixed-point sublattices of the positive-definite Leech lattice $\mathbb L$ \cite{HM16}.  We write $\mathbb L(-1)$ for its negative-definite version.

We use the following sign-reversed form of the Leech embedding lemma of Gaberdiel--Hohenegger--Volpato \cite[Appendix B.2]{GHV12} (the gap in the original argument concerning the primitivity of this embedding was fixed in ranks at most $21$ by Marquand--Muller \cite[Corollary 4.19]{MM25} and in full generality by Zheng \cite[Lemma 1.1]{zheng2025leechlemma}).

\begin{lemma} \label{lem:leech-embedding}
Let $(G,S)$ be a Leech pair.  If $\rk S+\ell(A_S)\leq24$, then $S$ admits a primitive $G$-equivariant embedding into $\mathbb L(-1)$, and the extended action of $G$ fixes the orthogonal complement of $S$ pointwise.
\end{lemma}

We also show the geometric meaning of saturation.

\begin{lemma}
\label{lem:saturation-on-surface}
Let a finite group $G$ act faithfully and symplectically on a
supersingular K3 surface $X$ in characteristic $7$.  Then the saturation of
$(G,S_G(X))$ acts symplectically on $X$ and extends the given
$G$-action.
\end{lemma}

\begin{proof}
Put $S=S_G(X)$ and $K=\NS(X)^G$; these are primitive orthogonal
complements in $\NS(X)$.  Let $\Gamma\subset A_S\oplus A_K$ be the isotropic subgroup defining the gluing from $S\oplus K$ to $\NS(X)$.  Every element of the saturation $\widetilde G$ acts trivially on $A_S$, so its direct sum with the
identity on $K$ fixes $\Gamma$ pointwise.  It therefore extends to an isometry of $\NS(X)$ acting trivially on $A_{\NS(X)}$.

The extension fixes a $G$-invariant ample class in $K$, hence preserves the ample cone, and its trivial action on $A_{\NS(X)}$ fixes the characteristic subspace.  Crystalline Torelli \cite[Theorem II]{Ogu83} and Jang's discriminant criterion \cite[Proposition 3.1]{Jan16} realize these extensions as a symplectic action of $\widetilde G$ on $X$ which extends the original $G$-action.
\end{proof}

We shall repeatedly use the following module estimate, which we state
for an arbitrary prime $p$ and apply only with $p=7$.

\begin{lemma}
\label{lem:cyclic-seven-modular-bound}
Let $p$ be a prime.
\begin{enumerate}
\item Let $C_p$ act on a free abelian group $M$
of finite rank.  If $M^{C_p}=0$, then
\[
\dim_{\BF_p}(M/pM)^{C_p}=\frac{\rk M}{p-1}.
\]
\item If $C_p$ acts on an integral lattice $L$ with $L^{C_p}=0$ and
acts trivially on $A_L$, then $A_L$ is annihilated by $p$ and
\[
\ell(A_L)\leq\frac{\rk L}{p-1}.
\]
\item Suppose that a finite group $G$ acts on an integral
lattice $L$, that $L^G=0$, and that $G$ acts trivially on $A_L$.  If
$C_p\leq G$ is a Sylow $p$-subgroup, then
\[
\ell_p(A_L)\leq\frac{\rk L-\rk L^{C_p}}{p-1}.
\]
\end{enumerate}
\end{lemma}

\begin{proof}
\textup{(1)} Let $g$ generate $C_p$.  Since $M^{C_p}=0$, the rational representation
$M\otimes\BQ$ has no trivial summand.  Hence, for some $d\geq0$,
\[
\det(T\operatorname{id}-g\mid M\otimes\BQ)=\Phi_p(T)^d,
\qquad \rk M=(p-1)d.
\]
It follows that $\abs{\det(g-1)}=\Phi_p(1)^d=p^d$.  Write
$\Phi_p(T)=p+(T-1)Q(T)$ with $Q(T)\in\BZ[T]$.  Since
$\Phi_p(g)=0$ on $M$, we have
\[
pM\subseteq(g-1)M.
\]
Thus $M/(g-1)M$ is annihilated by $p$ and has order $p^d$.
The Smith normal form of $g-1$ therefore has exactly $d$ diagonal
entries equal to $p$, with all remaining entries equal to $1$.
Since $g$ generates $C_p$, we have
\[
(M/pM)^{C_p}
=\{m+pM\mid m\in M,\ (g-1)m\in pM\}.
\]
The Smith normal form above then gives
\[
\dim_{\BF_p}(M/pM)^{C_p}=d=\frac{\rk M}{p-1}.
\]

\textup{(2)} Put $d=\rk L/(p-1)$.
For $x\in L^\vee$, triviality on $A_L$ gives $hx-x\in L$ for every
$h\in C_p$.  Since $L^{C_p}=0$, one has
$\sum_{h\in C_p}hx=0$ in $L\otimes\BQ$, and hence
$px=\sum_{h\in C_p}(x-hx)\in L$.
Thus $A_L$ is annihilated by $p$, and we have
\[
A_L=L^\vee/L\cong pL^\vee/pL\hookrightarrow(L/pL)^{C_p}.
\]
Hence $\dim_{\BF_p}A_L\leq d$.

\textup{(3)} Put $S=S_{C_p}(L)$ and let $g$ generate $C_p$.
Then $S^{C_p}=0$ and $(g-1)L^\vee\subseteq S$, so there is a map
\[
A_L[p]\longrightarrow S/(g-1)L,
\qquad x+L\longmapsto(g-1)x+(g-1)L.
\]
We claim that it is injective.  A class in its kernel has a representative
$x\in(L\otimes\BQ)^{C_p}$, after subtracting an element of $L$.
Let $t_1,\ldots,t_m$ represent the left cosets of $C_p$ in $G$.
Then $\sum_i t_i x$ is $G$-invariant, hence zero, so
\[
mx=\sum_i(x-t_i x)\in L.
\]
Since also $px\in L$ and $p\nmid m=[G:C_p]$, we obtain $x\in L$.
By the Smith normal form computation in (1),
\[
\abs{S/(g-1)S}=p^{\rk S/(p-1)}.
\]
Since $S/(g-1)L$ is a quotient of $S/(g-1)S$, the injection gives
\[
\ell_p(A_L)\leq\frac{\rk S}{p-1}
=\frac{\rk L-\rk L^{C_p}}{p-1}.
\]
\end{proof}

\section{Wild Groups on Finite-Height K3 Surfaces}
\label{sec:wild-degree-seven}
\label{sec:finite-height-seven}

We first deal with the finite-height case using the two non-unit
slopes of crystalline cohomology.

\begin{proposition}
\label{prop:finite-height-seven}
Let $X/k$ be a K3 surface of finite height over an algebraically closed
field $k$ of characteristic $7$, and let
$G\leq\Aut^s(X)$ be finite with $7\mid\abs G$.  Then the following
hold.
\begin{enumerate}
\item The height of $X$ is $6$, and every Sylow
$7$-subgroup of $G$ is cyclic of order $7$.
\item If $g\in G$ has order $7$, then there is a $g$-invariant
genus-one fibration for which the Kodaira type of a singular fibre,
the Picard number, and the action on the N\'eron--Severi group form
one of the following rows:
\[
\begin{array}{c|c|c}
\text{Kodaira type}&\rho(X)&
\det(T\operatorname{id}-g\mid\NS(X)\otimes\BQ)\\
\hline
II^*&10&(T-1)^{10}\\
III&4&(T-1)^4\\
III&10&(T-1)^4\Phi_7(T).
\end{array}
\]
\item One has $G\cong C_7$ or $C_{14}$, and both groups occur.
If $G\cong C_{14}$, then its order-$7$ subgroup occurs in the last row
of part (2), namely the row with a singular fibre of type $III$ and
$\rho(X)=10$.
More precisely, the minimal resolution of
\begin{equation}
\label{eq:finite-height-c14}
y^2=x^3+s^4x^2+(t^7s-ts^7)x
\subset\mathbb P(1,1,4,6)
\end{equation}
has height $6$ and carries a symplectic $C_{14}$-action.
\end{enumerate}
\end{proposition}

\begin{proof}
Choose an element $g\in G$ of order $7$.  By
\cite[Proposition 8.1(1),(2)]{Keu16}, the quotient
$X/\langle g\rangle$ is rational and $X$ admits a
$g$-invariant genus-one fibration
\[
\psi\colon X\longrightarrow\mathbb P^1
\]
with a distinguished $g$-invariant fibre $F_0$.  By \cite[Proposition 8.1(3)--(5)]{Keu16}, the fibre $F_0$ is of type
$II^*$ or $III$, and the corresponding characteristic polynomial on $H^2_{\mathrm{\acute et}}(X,\BQ_\ell)$ is respectively
\begin{equation}
\label{eq:keum-order-seven-characters}
\det(T\operatorname{id}-g\mid H^2_{\mathrm{\acute et}}(X,\BQ_\ell))
=(T-1)^{10}\Phi_7(T)^2
\quad\text{or}\quad
(T-1)^4\Phi_7(T)^3.
\end{equation}
Let
\[
T_\ell(X)=(\NS(X)\otimes\BQ_\ell)^\perp
 \subset H^2_{\mathrm{\acute et}}(X,\BQ_\ell(1))
\]
be the $\ell$-adic transcendental space.  By Keum \cite[Proposition 2.2(4)]{Keu16}, $T_\ell(X)^g=0$ and hence
\[
\dim H^2_{\mathrm{\acute et}}(X,\BQ_\ell)^g
=\rk\NS(X)^g.
\]

Let $W=W(k)$, $K_0=\operatorname{Frac}W$, and let $D$ be the part of $H^2_{\mathrm{cris}}(X/W)\otimes K_0$ with slope $<1$.  If $h$ is the formal Brauer height, then $\dim_{K_0}D=h$, while all divisor classes have slope $1$; see \cite[\S 1]{Jan17}.  Thus $D$ lies in the crystalline transcendental space
\[
T_{\mathrm{cris}}(X)
=(\NS(X)\otimes W)^\perp\otimes_W K_0.
\]
The characteristic polynomial on total cohomology is independent of the choice between $\ell$-adic and crystalline cohomology \cite[Proposition 2.1(1)]{Keu16}.  Removing the common N\'eron--Severi factor shows that $T_{\mathrm{cris}}(X)$ and $T_\ell(X)$ have the same characteristic polynomial.  Hence $D^g=0$.

Since the polynomial $\Phi_7(T+1)$ is Eisenstein over $W$, $\Phi_7$ is irreducible over $K_0$ of degree 6. It follows that $6\mid h$.  Since $1\le h\le 10$, one has $h=6$.

Let $P\in\Syl_7(G)$.  Every nonidentity element of $P$ has order
$7$ by \cite[Theorem 8.2]{Keu16}, and the preceding argument gives
characteristic polynomial $\Phi_7$ on $D$, hence trace $-1$.
Since $D^P=0$, character averaging gives
\[
0=\dim D^P=\frac{6-(\abs P-1)}{\abs P}
=\frac7{\abs P}-1.
\]
Thus $P\cong C_7$, proving part (1).

For part (2), \eqref{eq:keum-order-seven-characters} gives
$\rk\NS(X)^g=10$ or $4$ according as the fibre $F_0$ chosen above
has type $II^*$ or $III$.  Since $\rho(X)\leq22-2h=10$ and every
nontrivial rational $C_7$-summand of $\NS(X)\otimes\BQ$ has
dimension $6$, these give exactly the three rows in the statement.

For a tame subgroup $H\leq G$, \cite[Proposition 1.3]{WZroot}
gives $\rk S_H(X)\leq9$ and a primitive embedding into the K3
lattice, hence $\rk S_H(X)+\ell(A_{S_H(X)})\leq22$.
The pair $(H,S_H(X))$ is a Leech pair, so Lemma
\ref{lem:leech-embedding} and \cite[Table 1]{HM16} give
$\abs H\leq2$.  Applying this to the Sylow subgroups away from
$7$ gives $\abs G=7$ or $14$.
In the dihedral case, choose an involution $a$ with
$aga^{-1}=g^{-1}$.  Its N\'eron--Severi preserving equivariant
lifting \cite[Corollary 3.4]{Jan17} is symplectic on the generic
fibre, so $a$ acts trivially on the transcendental lattice.
Comparison gives $a|_D=1$, contradicting the inversion relation
because $g|_D$ has order $7$.  Thus $G\cong C_7$ or $C_{14}$.

Suppose that $G\cong C_{14}$, and let $a$ be its involution.  The coinvariant lattice of the tame symplectic action of $a$ is H\"ohn--Mason no.~2 and has rank $8$ \cite[Table 1]{HM16}.  Since this lattice lies in $\NS(X)$, one has $\rho(X)\geq8$, so the row with $\rho(X)=4$ is impossible.  If the specific $g$-invariant fibre were of type $II^*$, then $X^g$ would be its support: the alternative in which $X^g$ is a point occurs only for type $III$ \cite[Proposition 8.1(2),(4)]{Keu16}.  Since $a$ commutes with $g$, it permutes the components of this fibre.  Their dual graph is the affine
$E_8$ Dynkin diagram, which has no nontrivial automorphism, so $a$ preserves every component.  On the trivalent component it fixes the three intersection points and therefore acts identically.  This contradicts the fact that a tame symplectic involution has only isolated fixed points \cite[Theorem 3.3]{DK09}.  Thus only the row with a singular fibre of type $III$ and $\rho(X)=10$ can occur.

It remains to verify occurrence.  Equation \eqref{eq:finite-height-c14}, which is among the elliptic models of
\cite[Examples 5.8, case 7-(5)]{DK01}, has affine form
\[
y^2=x^3+x^2+qx,\qquad q=t^7-t.
\]
Its invariants are $\Delta=16q^2(1-4q)$ and $c_4=16-48q$.  Thus its minimal resolution
has seven fibres of type $I_2$, seven of type $I_1$, and a fibre of
type $III$ at infinity.  The relative Weierstrass model has the eight
$A_1$ singularities belonging to the seven $I_2$ fibres and the
$III$ fibre.  The coarse weighted hypersurface has one additional
$A_1$ point at its intersection with the ambient singular line
$\{s=t=0\}$.  The Euler numbers add to $24$,
and the global Weierstrass coefficients have degrees $4$ and $8$, so
the resolution is an elliptic K3 surface.

The transformations
\[
g(t,x,y)=(t+1,x,y)
\quad\text{and}\quad
\tau(x,y)=\left(\frac qx,-\frac{qy}{x^2}\right)
\]
are respectively the base translation and translation by the
$2$-torsion section $(0,0)$.  They extend to the minimal elliptic
surface, commute, and preserve $dt\wedge dx/y$.  Hence they generate
a symplectic group $C_7\times C_2\cong C_{14}$.

The following lemma computes its height.
\end{proof}
\begin{lemma}
Over $\BF_7$, the minimal resolution of the hypersurface
$\bar X=\{f=0\}\subset\mathbb P(1,1,4,6)$, where
\[
f=y^2-x^3-s^4x^2-(t^7s-ts^7)x,
\]
is a K3 surface of height $6$.
\end{lemma}
\begin{proof}
Let $\pi:X\to\bar X$ be the minimal resolution.  The hypersurface
$\bar X$ has only rational double points of type $A_1$, so $\pi$ is crepant:
\[
K_X=\pi^*K_{\bar X}.
\]
As shown above, $X$ is a K3 surface.  Over an algebraic closure,
Yobuko \cite{Yob19} identifies its height $h$ with its quasi-$F$-split
height.  Crepant birational invariance
\cite[Proposition 3.46, Corollary 3.47]{KTT22} identifies this with
that of $\bar X$.

Set $p=7$ and
\[
S=\mathbb F_7[s,t,x,y],
\qquad
\deg(s)=\deg(t)=1,\quad
\deg(x)=4,\quad
\deg(y)=6,
\]
and put
\[
\mathfrak m=(s,t,x,y),
\qquad
\mathfrak m^{[7]}=(s^7,t^7,x^7,y^7).
\]
The weighted projective space $\mathbb P(1,1,4,6)$ is well formed
\cite[Definition 3.1]{KTY25II}, since any three weights have greatest
common divisor one.  Moreover, $\deg f=1+1+4+6=12$, and $\bar X$
meets the ambient singular locus only in codimension two.  Thus the
weighted projective criterion \cite[Theorem 3.16(1)]{KTY25II} applies.

Let $u:F_*S\rightarrow S$ be the Frobenius trace: \[ u(F_*(s^{a}t^{b}x^{c} y^{d}))=\begin{cases}
s^{(a-6)/7}t^{(b-6)/7}x^{(c-6)/7}y^{(d-6)/7},
&
a,b,c,d\equiv6\pmod7,\\[2mm]
0,&\text{otherwise}.
\end{cases}\]
For $a\in S$, write $a=\sum_i M_i$ as a sum of its nonzero
monomial terms, including their coefficients.  Lift each $M_i$ to
$\widetilde M_i\in\BZ[s,t,x,y]$ by choosing its coefficient in
$\{0,\ldots,6\}$, and put $\widetilde a=\sum_i\widetilde M_i$.
Define
\[
\Delta_1(a)=
\left(\frac{\widetilde a^{\,p}-\sum_i\widetilde M_i^{\,p}}{p}\right)
\bmod p\in S,
\]
where the numerator is divisible by $p$ in $\BZ[s,t,x,y]$ and the
division is performed before reduction modulo $p$.  Recursively set
\[
H_1=f^{p-1},
\qquad
H_{m+1}=u\bigl(F_*(\Delta_1(f^{p-1})H_m)\bigr).
\]
By \cite[Theorem 5.8]{KTY25I} and
\cite[Theorem 3.16(1)]{KTY25II}, it suffices to verify
\[
H_m\in \mathfrak m^{[7]}, \qquad 1\leq m\leq5, \qquad H_6\notin \mathfrak m^{[7]}.
\]

We have $\deg \Delta_1(f^{p-1})=504$ and $H_1\in S_{72}$. The trace sends the weighted-degree-$576$ part of $S$ to $S_{72}$, so induction gives
\[
H_m\in S_{72}
\qquad\text{for every }m\geq1.
\]
The only monomial of weighted degree $72$ outside $\mathfrak m^{[7]}$ is $s^6t^6x^6y^6$. We denote the coefficient of $s^6t^6x^6y^6$ in $H_m$ by $c_m$. The accompanying script
\path{anc/height_six_verification.py} computes the recurrence using
exact arithmetic and gives
\[
\begin{array}{c|rrrrrr}
m&1&2&3&4&5&6\\ \hline
\#\operatorname{supp}(H_m)
&210&209&337&455&536&587\\
c_m&0&0&0&0&0&6.
\end{array}
\]
Thus
\[
H_m\in\mathfrak m^{[7]}
\quad(1\leq m\leq5),
\qquad
H_6\notin\mathfrak m^{[7]}.
\]
Hence $\bar X$ has quasi-$F$-split height $6$, unchanged after extension
to an algebraic closure \cite[Theorem 5.13]{KTY25I}.  The crepant
comparison and Yobuko's theorem then give $h(X)=6$.
\end{proof}
\section{Finite Symplectic Groups on Supersingular K3 Surfaces}
\label{sec:fine-classification}

\subsection{Wild automorphisms}
\label{sec:geometric-bound}

From this point on we use the realization criterion
\cite[Proposition 1.1]{WZroot}: a Leech pair $(G,S)$ is realized on a
supersingular K3 surface of Artin invariant $\sigma$ if and only if
$S$ admits a primitive embedding into $N_{7,\sigma}$; whenever this
holds, the pair is realized on every surface in that Artin stratum.

\Needspace{5\baselineskip}
\begin{theorem}
\label{thm:wild-degree-seven}
In characteristic $7$, the following hold.
\begin{enumerate}
\item Let $X$ be a supersingular K3 surface, and let
$g\in\Aut^s(X)$ have order $7$.  Then, with $T$ an indeterminate,
\[
\det(T\operatorname{id}-g\mid\NS(X)\otimes\BQ)
=(T-1)^4\Phi_7(T)^3.
\]
The saturation of
$(\langle g\rangle,S_{\langle g\rangle}(X))$ is H\"ohn--Mason no.~52,
whose group is $F_{21}=C_7\rtimes C_3$.  In particular, the given
$\langle g\rangle$-action extends on $X$ to a symplectic
$F_{21}$-action.
\item If a finite group $G$ acts faithfully and symplectically on a
supersingular K3 surface, then $\nu_7(G)\leq1$, where
$7^{\nu_7(G)}\mathrel\Vert\abs G$.
\item If $X$ is supersingular of Artin invariant $\sigma\geq4$, then
every finite-order automorphism of $X$ is tame.
\end{enumerate}
\end{theorem}

\begin{proof}[Proof of Theorem \ref{thm:wild-degree-seven}(1)]
On $\NS(X)\otimes\BQ$, the rational characteristic polynomial has the
form
\[
\det(T\operatorname{id}-g\mid\NS(X)\otimes\BQ)
=(T-1)^{22-6m}\Phi_7(T)^m,
\qquad 1\leq m\leq3.
\]
The coinvariant lattice $S_{\langle g\rangle}(X)$ is negative-definite,
even, and rootless of rank $6m$.  The action of $\langle g\rangle$
on it is faithful and fixed-point-free and is trivial on its
discriminant group.  Thus
$(\langle g\rangle,S_{\langle g\rangle}(X))$ is a Leech pair.  Lemma
\ref{lem:cyclic-seven-modular-bound} shows that its discriminant group
is annihilated by $7$ and has length at most $m$.  Hence
\[
\rk S_{\langle g\rangle}(X)
 +\ell(A_{S_{\langle g\rangle}(X)})
\leq7m\leq21.
\]
Lemma \ref{lem:leech-embedding} embeds this pair primitively into
$\mathbb L(-1)$, and the extended action has
$S_{\langle g\rangle}(X)$ as its coinvariant lattice.  Every
order-$7$ element of $\Co_0$ has coinvariant rank at least $18$
\cite[Table 1]{HM16}.  Hence $m=3$.

Every element of the saturation extends by the identity on the
orthogonal complement, so its fixed lattice in $\mathbb L(-1)$ has
rank $24-18=6$.  The unique
saturated H\"ohn--Mason entry of fixed rank $6$ whose order is
divisible by $7$ is no.~52, with group $F_{21}$
\cite[Table 1]{HM16}.

Lemma \ref{lem:saturation-on-surface} now realizes this saturation as
a symplectic $F_{21}$-action on $X$ containing the given
automorphism.
\end{proof}

\begin{proof}[Proof of Theorem \ref{thm:wild-degree-seven}(2) and (3)]
Let $G\leq\Aut^s(X)$ be finite and let $P\in\Syl_7(G)$ be nontrivial.
An element of order $49$ would force a $\Phi_{49}$-summand of
dimension $\varphi(49)=42$ in the faithful rational representation on
$\NS(X)$, which has rank $22$.  Hence $P$ has exponent $7$, and every nonidentity element of
$P$ has trace $1$ on $\NS(X)\otimes\BQ$ by part (1).  Character averaging
gives
\[
\rk\NS(X)^P
=\frac{22+(\abs P-1)}{\abs P}
=1+\frac{21}{\abs P}.
\]
This is an integer, so the nontrivial $7$-group $P$ has order $7$.
This proves part (2).

Let $g$ generate $P$.  The characteristic polynomial in part (1)
shows that $\NS(X)^g$ has rank $4$, while Lemma
\ref{lem:cyclic-seven-modular-bound} gives
$\ell_7(A_{S_{\langle g\rangle}(X)})\leq3$.  Since
$S_{\langle g\rangle}(X)$ and $\NS(X)^g$ are primitive orthogonal
complements, discriminant gluing gives
\[
2\sigma=\ell_7(A_{\NS(X)})
\leq\ell_7(A_{S_{\langle g\rangle}(X)})
 +\ell_7(A_{\NS(X)^g})\leq3+4=7.
\]
Thus $\sigma\leq3$ whenever a wild finite
symplectic group acts.

Finally, let $h\in\Aut(X)$ have finite order $n$ divisible by $7$.
Then $h^{n/7}$ has order $7$ and acts trivially on
$H^0(X,\Omega_X^2)$, because a field of characteristic $7$ has no
nontrivial seventh root of unity.  Thus it is a wild symplectic
automorphism, which is impossible when $\sigma\geq4$.  This proves
part (3).
\end{proof}

\subsection{Wild groups}

We now retain the group, rather than only its $7$-part.  We first
classify the wild actions through their saturated pairs, then record
the corresponding tame higher-Artin lists.

\begin{proposition}
\label{prop:uniform-leech-seven}
Let $X$ have Artin invariant $\sigma\geq2$, let
$G\leq\Aut^s(X)$ be finite with $7\mid\abs G$, and put
$r=\rk\NS(X)^G$.  Then
\[
2\leq r\leq4,
\qquad
\ell_7(A_{S_G(X)})\leq3,
\qquad
\ell_q(A_{S_G(X)})\leq r\quad(q\neq7).
\]
Consequently $(G,S_G(X))$ embeds into $\mathbb L(-1)$, and the
possible quadruples $(\sigma,r,i,\widetilde G)$, where $i$ is the
H\"ohn--Mason number of $S_G(X)$, are contained in the following
table:
\[
\begin{array}{c|c|c|c}
\sigma&r&\text{H\"ohn--Mason number}&\widetilde G\\
\hline
2,3&4&52&F_{21}=C_7\rtimes C_3\\
2&3&77&L_2(7)\\
2&2&129&C_2\times L_2(7).
\end{array}
\]
\end{proposition}

\begin{proof}
Let $P\in\Syl_7(G)$.  Theorem \ref{thm:wild-degree-seven}(2) gives
$P\cong C_7$, and parts (1) and (3) give
\[
\sigma\in\{2,3\},
\qquad
\rk\NS(X)^P=4.
\]
A $G$-invariant ample class and the inclusion
$\NS(X)^G\subset\NS(X)^P$ give $1\leq r\leq4$.  Rational
orthogonal decomposition gives
\[
\rk S_G(X)^P=\rk\NS(X)^P-\rk\NS(X)^G=4-r.
\]
Although this fixed sublattice need not vanish, the final assertion of
Lemma \ref{lem:cyclic-seven-modular-bound} applies to $S_G(X)$ and
gives
\[
\ell_7(A_{S_G(X)})
\leq\frac{\rk S_G(X)-\rk S_G(X)^P}{6}
=\frac{(22-r)-(4-r)}6=3.
\]

For $q\neq7$, the lattice $\NS(X)\otimes\BZ_q$ is unimodular.
Primitive orthogonal gluing with the rank-$r$ lattice $\NS(X)^G$
therefore gives $\ell_q(A_{S_G(X)})\leq r$.  Hence
\[
\rk S_G(X)+\ell(A_{S_G(X)})
\leq22-r+\max\{3,r\}\leq24,
\]
so Lemma \ref{lem:leech-embedding} embeds $(G,S_G(X))$ into
$\mathbb L(-1)$.

Primitive discriminant gluing also gives
\[
2\sigma=\ell_7(A_{\NS(X)})
\leq\ell_7(A_{S_G(X)})+\ell_7(A_{\NS(X)^G})
\leq\ell_7(A_{S_G(X)})+r.
\]
The saturation has fixed rank $24-\rk S_G(X)=r+2$ in
$\mathbb L(-1)$.  Among the entries of this fixed rank whose group
order is divisible by $7$, \cite[Table 1]{HM16} and the displayed
inequality leave precisely the candidate rows above; in particular,
$r\geq2$.
\end{proof}

\begin{remark}
Wang--Zheng's root-pair method studies the roots that may appear after
an equivariant extension from $N_{p,\sigma}$ to $N_{p,1}$
\cite[Proposition 1.6 and \S\S 3.5, 6.2]{WZroot}.  The preceding proof
does not make this extension.  The bound $\nu_7(G)\leq1$ is obtained
directly on $N_{7,\sigma}$ in Theorem
\ref{thm:wild-degree-seven}(2), and Lemma
\ref{lem:leech-embedding} is applied to the original rootless
coinvariant lattice $S_G(N_{7,\sigma})$.  Thus the argument reaches
H\"ohn--Mason nos.~52, 77, and 129 without using a root pair.  This
does not mean that $S_G(N_{7,1})$ is always rootless: its root lattice
can depend on the chosen equivariant extension.
\end{remark}

We next prove that the three remaining saturated pairs have the
primitive embeddings needed for sharpness.  Let $S_i$ be the coinvariant
lattice of H\"ohn--Mason no.~$i$.  Here and below, we use the Conway--Sloane
notation for finite quadratic forms \cite[Chapter 15, \S 7]{CS99}; in
particular, $7^{\varepsilon m}$ is a $7$-elementary form of length $m$,
with $\varepsilon\in\{+,-\}$ recording its determinant square class.
The three relevant saturated pairs have the following data
\cite[Table 1]{HM16}:
\[
\begin{array}{c|c|c|c}
i&\rk S_i&\widetilde G_i&q_{S_i}\\
\hline
52&18&F_{21}&7^{+3}\\
77&19&L_2(7)&4_1^{+1}\oplus7^{+2}\\
129&20&C_2\times L_2(7)&2_{\mathrm{II}}^{+2}\oplus7^{+2}.
\end{array}
\]

\begin{proposition}
\label{prop:seven-primitive-embeddings}
There are primitive embeddings
\[
S_{52}\hookrightarrow N_{7,2},
\qquad
S_{52}\hookrightarrow N_{7,3},
\qquad
S_{77}\hookrightarrow N_{7,2},
\qquad
S_{129}\hookrightarrow N_{7,2}.
\]
\end{proposition}

\begin{proof}
Let $U$ be the hyperbolic plane and let $B_7$ have Gram matrix
\[
\begin{pmatrix}
-2&-1\\
-1&-4
\end{pmatrix}.
\]
Choose $K$ of signature $(1,21-\rk S_i)$ from the following table:
\[
\begin{array}{c|c|c|c|c}
i&\sigma&K&(q_K)_2&(q_K)_7\\
\hline
52&2&U\oplus B_7&0&7^{-1}\\
52&3&U(7)\oplus B_7&0&7^{+3}\\
77&2&U(7)\oplus\langle-4\rangle&-4_1^{+1}&7^{-2}\\
129&2&U(14)&2_{\mathrm{II}}^{+2}&7^{-2}.
\end{array}
\]
In each row, $(q_K)_2\cong-(q_{S_i})_2$.  Glue $S_i\oplus K$
along the graph of an anti-isometry between these parts
\cite[Proposition 1.4.1]{Nik79}.  The resulting even lattice $N$
contains both factors primitively, since the graph meets each
discriminant summand trivially; for $i=52$, this is simply the direct sum.  It has
signature $(1,21)$ and discriminant form
\[
q_N\cong(q_{S_i})_7\oplus(q_K)_7\cong
\begin{cases}
7^{-4},&\sigma=2,\\
7^{+6},&\sigma=3.
\end{cases}
\]
These are precisely $q_{N_{7,\sigma}}$
\cite[Proposition 2.14]{WZroot}.  The Rudakov--Shafarevich uniqueness
theorem \cite[\S 1]{RS81} gives $N\cong N_{7,\sigma}$, proving the
required primitive embeddings.
\end{proof}

We now pass from the saturated pairs and their embeddings to the wild
groups themselves.

\begin{theorem}
\label{thm:finite-groups-seven}
Let a finite wild group $G$ act faithfully and symplectically on a
supersingular K3 surface $X$ of Artin invariant $\sigma\geq2$ in
characteristic $7$.  Then the triple consisting of $G$, the fixed
rank, and the saturated pair occurs in one of the following rows:
\[
\begin{array}{c|c|c|c}
\sigma&G&\rk\NS(X)^G&(\text{H\"ohn--Mason no.},\widetilde G)\\
\hline
2,3&C_7,\ F_{21}&4&(52,F_{21})\\
2&L_2(7)&3&(77,L_2(7))\\
2&C_{14},\ C_2\times F_{21},\ C_2\times L_2(7)
&2&(129,C_2\times L_2(7)).
\end{array}
\]
Every such $G$-action extends on the same surface to its saturation
$\widetilde G$.  Conversely, each pair $(\sigma,G)$ displayed in the
table is realized on every supersingular K3 surface of Artin invariant
$\sigma$, with the fixed rank and saturated pair shown.
\end{theorem}

\begin{proof}
By Proposition \ref{prop:uniform-leech-seven}, $G$ embeds in
$C_2\times L_2(7)$.  The subgroup structure of $L_2(7)$
\cite{Atlas85} gives projection $H\cong C_7,F_{21}$, or $L_2(7)$;
since $\operatorname{Hom}(H,C_2)=0$, the group $G$ is isomorphic to
$H$ or $C_2\times H$.

Theorem \ref{thm:wild-degree-seven}(1) gives saturation no.~52 for
$C_7$.  For $F_{21}$ and $L_2(7)$, character averaging with the
element counts in \cite{Atlas85} and the traces in Theorem
\ref{thm:wild-degree-seven}(1) and \cite[Theorem 3.3]{DK09} gives
\[
\begin{aligned}
\rk\NS(X)^{F_{21}}
&=\frac{22+14\cdot4+6\cdot1}{21}=4,\\
\rk\NS(X)^{L_2(7)}
&=\frac{22+21\cdot6+56\cdot4+42\cdot2+48\cdot1}{168}=3.
\end{aligned}
\]
Proposition \ref{prop:uniform-leech-seven} therefore identifies their
saturations as nos.~52 and 77, respectively.  Each remaining group
contains an element of order $14$, so its saturation is no.~129,
since neither $F_{21}$ nor $L_2(7)$ contains such an element.
This gives the fixed ranks and saturated pairs in the table; the
extension on the same surface is Lemma
\ref{lem:saturation-on-surface}.

Conversely, Proposition \ref{prop:seven-primitive-embeddings} and
\cite[Proposition 1.1]{WZroot} realize the three saturated pairs on
every surface in the indicated Artin strata.  Restriction to the
subgroups in each corresponding row preserves the fixed rank,
hence the coinvariant lattice, and realizes every listed case.
For $\sigma=2$, existence already follows from
\cite[Remark 9.2]{OS24}.
\end{proof}

\begin{remark}
Inside $2^4{:}A_7\leq M_{23}$, the stabilizer $L_2(7)\leq A_7$
of a nonzero vector in $2^4$ centralizes that vector, giving
$C_2\times L_2(7)\hookrightarrow M_{23}$ \cite{Atlas85}.
\end{remark}

The classification also has the following cohomological consequence.

\begin{corollary}
Let $X$ be a supersingular K3 surface of Artin invariant
$\sigma\geq2$ in characteristic $7$, and let
$G\leq\Aut^s(X)$ be finite.  Then there is an embedding
$G\hookrightarrow M_{23}$ under which, for every prime $\ell\neq7$,
the character on $H^*_{\mathrm{\acute et}}(X,\BQ_\ell)$ is the
restriction of the natural twenty-four-point permutation character
of $M_{23}$.
\end{corollary}

\begin{proof}
For tame $G$, both assertions are due to Dolgachev--Keum
\cite[Theorems 3.3, 4.7 and Proposition 4.1]{DK09}.
For wild $G$, Theorem \ref{thm:finite-groups-seven} and the
preceding remark give the embedding into $M_{23}$ and the possible
element orders $1,2,3,4,6,7,14$.  The tame traces are known,
and the order-$7$ trace is $2+(4-3)=3$ by Theorem
\ref{thm:wild-degree-seven}(1).

For an element $g$ of order $14$, Theorem
\ref{thm:finite-groups-seven} gives
\[
\dim H^*_{\mathrm{\acute et}}(X,\BQ_\ell)^{\langle g\rangle}
=2+\rk\NS(X)^{\langle g\rangle}=4.
\]
Since $X$ is supersingular, this representation is the scalar extension of
$\BQ^2\oplus(\NS(X)\otimes\BQ)$, so the six generators of
$\langle g\rangle$ have the same trace $\tau$.
Character averaging gives
\[
14\cdot4=24+8+6\cdot3+6\tau,
\qquad \tau=1.
\]
Thus the traces are
\[
\begin{array}{c|ccccccc}
\operatorname{ord}(g)&1&2&3&4&6&7&14\\
\hline
\operatorname{Tr}(g\mid H^*_{\mathrm{\acute et}}(X,\BQ_\ell))&24&8&6&4&2&3&1,
\end{array}
\]
which agree with the natural twenty-four-point permutation
character of $M_{23}$ \cite{Muk88}.
\end{proof}

\subsection{Tame groups by Artin invariant}

The superspecial table used in the tame reduction is
\cite[Appendix A]{WZroot}; see also
\cite[Propositions 3.9 and 3.10]{WZroot}.

\begin{lemma}
\label{lem:tame-seven-embedding}
Let $S$ be the coinvariant lattice of a saturated H\"ohn--Mason pair in
the characteristic-$7$ superspecial list of
\cite[Appendix A]{WZroot}, suppose that
$7\nmid\abs{A_S}$, and put $c=22-\rk S$.  For $2\leq\sigma\leq10$,
the lattice $S$ admits a primitive embedding into $N_{7,\sigma}$ if
and only if either $c>2\sigma$, or $c=2\sigma$ and
\[
\left(\frac{\abs{A_S}}7\right)=(-1)^\sigma.
\]
\end{lemma}

\begin{proof}
Fix a primitive embedding $S\hookrightarrow N_{7,1}$ with complement
$K_1$.  Primitive discriminant gluing and uniqueness of supersingular
K3 lattices show that $S\hookrightarrow N_{7,\sigma}$ exists if and
only if there is an even lattice $K_\sigma$ of signature $(1,c-1)$ with
\[
q_{K_\sigma}\cong-q_S\oplus q_{N_{7,\sigma}}.
\]
In particular, $c\geq\ell_7(A_{K_\sigma})=2\sigma$ is necessary.

The local conditions away from $7$ are already satisfied by $K_1$,
since the determinant changes by the square $7^{2\sigma-2}$.  The
signature condition is also unchanged, since $q_{N_{7,\sigma}}$ has
the same Brown invariant for every $\sigma$.  Thus Nikulin's criterion
\cite[Theorem 1.10.1]{Nik79} gives existence when $c>2\sigma$.

When $c=2\sigma$, the remaining local condition compares the unit
part $\det K_\sigma/7^{2\sigma}=-\abs{A_S}$ with the determinant of
the $7$-primary discriminant form.  Since $q_{N_{7,\sigma}}$ is
$7^{+2\sigma}$ for odd $\sigma$ and $7^{-2\sigma}$ for even $\sigma$,
this condition is
\[
\left(\frac{-\abs{A_S}}7\right)=(-1)^{\sigma+1},
\qquad
\text{equivalently}\qquad
\left(\frac{\abs{A_S}}7\right)=(-1)^\sigma.
\]
Nikulin's criterion also proves sufficiency.
\end{proof}

\begin{theorem}
\label{thm:tame-groups-seven}
Let a finite tame group $G$ act faithfully and symplectically on a
supersingular K3 surface $X$ of Artin invariant $\sigma\geq2$ in
characteristic $7$.  Then the possible H\"ohn--Mason numbers of its
saturation are exactly
\[
\begin{array}{c|l}
\sigma&\text{H\"ohn--Mason numbers}\\
\hline
2&1,2,3,4,6,7,9,10,13,15,17,18,19,20,24,26,28,29,31,\\
 &37,43,45,46,55\\
3&1,2,3,4,6,7,9,10,13\\
4&1,2,3,4,6,9\\
5,6&1,2\\
7,8,9,10&1.
\end{array}
\]
Every listed saturated pair, and every tame subgroup of its saturated
group, acts symplectically on every surface in the corresponding
stratum.
\end{theorem}

\begin{proof}
For tame $G$, \cite[Proposition 1.2 and Lemma 4.1]{WZroot} give
$\rk S_G(X)\leq22-2\sigma$, $7\nmid\abs{A_{S_G(X)}}$, and a
saturated pair in the characteristic-$7$ superspecial list.
By Lemma \ref{lem:tame-seven-embedding}, all entries of rank
strictly below $22-2\sigma$ occur; at equality, one retains those
satisfying the Legendre condition.  Evaluating this criterion using
\cite[Table 1]{HM16} and the corrected discriminant forms in
\cite[Appendix A]{WZroot} gives exactly the displayed lists.
The realization criterion \cite[Proposition 1.1]{WZroot} realizes
each listed pair on every surface in its Artin stratum, and
restriction gives the assertion for subgroups.
\end{proof}

\begin{proof}[Proof of Theorem \ref{thm:all-wild-seven}]
If $X$ has finite height, Proposition \ref{prop:finite-height-seven}
gives the first case and realizes both groups on the surface
\eqref{eq:finite-height-c14}.  The case of Artin invariant one and
its converse are \cite[Theorem 1.2(4) and Corollary 1.3]{OS24}.
For Artin invariant at least two, Theorem
\ref{thm:finite-groups-seven} gives the remaining cases and their
realizations, while Theorem \ref{thm:wild-degree-seven}(3) excludes
Artin invariant at least four and proves the final assertion.
\end{proof}

\begin{proof}[Proof of Proposition \ref{prop:all-supersingular-seven}]
For $\sigma=1$, this is \cite[Theorem 1.2(4)]{OS24}.
For $\sigma\geq2$, the displayed list follows from Theorems
\ref{thm:all-wild-seven} and \ref{thm:tame-groups-seven}, and
realization on every surface in the corresponding Artin stratum
follows from Theorems \ref{thm:finite-groups-seven} and
\ref{thm:tame-groups-seven} by restriction to subgroups.
\end{proof}

\appendix

\section{An Explicit Supersingular Family}
\label{app:explicit-seven}

We give an elliptic $F_{21}$-family of generic Artin invariant two.
The special
members with $\mu=0,1$ already occur in \cite[Examples 5.8]{DK01}.

\begin{proposition}
\label{prop:explicit-seven-family}
Let $k$ be algebraically closed of characteristic $7$, put
$u=t^7-t$, and let $J_\mu$ be the relative minimal smooth model of the
Weierstrass surface
\[
J_\mu:\qquad y^2=x^3+ux+\mu,
\qquad \mu\in k.
\]
Then $J_\mu$ is a supersingular elliptic K3 surface.  The transformations
\[
g(x,y,t)=(x,y,t+1),
\qquad
h(x,y,t)=(4x,y,2t)
\]
extend to symplectic automorphisms of $J_\mu$ and satisfy
\[
g^7=h^3=1,
\qquad hgh^{-1}=g^2.
\]
They generate $F_{21}=C_7\rtimes C_3$, and $J_\mu$ has Artin
invariant two for $\mu$ in a nonempty Zariski-open subset of $k$.
\end{proposition}

\begin{proof}
The global Weierstrass equation is
\[
y^2=x^3+(t^7s-ts^7)x+\mu s^{12}
\subset\mathbb P(1,1,4,6).
\]
The coefficient of $x$ has only simple zeros on $\mathbb P^1$, so the
Weierstrass fibration is everywhere minimal with fundamental line
bundle $\mathcal O_{\mathbb P^1}(2)$.  Its minimal smooth model is
therefore a K3 surface.  The displayed coordinate changes preserve
the equation and the two-form
$dt\wedge dx/y$, and their relations are immediate.

The $F_{21}$-action excludes finite height by Proposition
\ref{prop:finite-height-seven}(3).  Thus $J_\mu$ is supersingular
\cite[Theorem 2.3]{Lie15}.  The surface $J_0$ is the
classical Artin-invariant-one model, isomorphic to the Fermat quartic
\cite[p.~311]{DK09}.

By Theorem \ref{thm:wild-degree-seven}(1), $\NS(J_\mu)^g$ has rank
$4$, and Lemma \ref{lem:cyclic-seven-modular-bound} gives
$\ell_7(A_{S_{\langle g\rangle}(J_\mu)})\leq3$.  The fibre and zero
section span a unimodular hyperbolic plane in $\NS(J_\mu)^g$, so
$\ell_7(A_{\NS(J_\mu)^g})\leq2$.  Discriminant gluing therefore gives
\[
2\sigma(J_\mu)\leq3+2=5,
\]
hence $\sigma(J_\mu)\leq2$.

For $\mu\neq0$, the $j$-map is
\[
j_\mu(t)=1728\,\frac{4u^3}{4u^3+27\mu^2},
\qquad
\Delta_\mu=-16(4u^3+27\mu^2).
\]
The twenty-one simple finite zeros of the discriminant give fibres of
type $I_1$, whereas the fibre at infinity is of type $III$.  An
equivalence of two marked fibrations therefore fixes infinity and
induces $t\mapsto at+b$.  If $J_\mu$ and $J_\nu$ are
equivalent as marked fibrations, equality of their $j$-maps gives
$u(at+b)=c\,u(t)$, and comparison of coefficients gives
\[
a^6=1,
\qquad b^7-b=0,
\qquad c=a,
\qquad \nu^2=a^3\mu^2.
\]
Only finitely many parameters define equivalent marked fibrations.
If the geometric generic K3 surface had Artin invariant one, its
uniqueness \cite[Corollary 7.14]{Ogu79} would trivialize the underlying
family after a finite base change and shrinking the base.  Since
$H^1(J_\mu,\mathcal O_{J_\mu})=0$, the Picard scheme is discrete, so the fibre
and zero-section classes would be constant.  The fibre linear system
and the unique effective zero-section representative then recover an
isotrivial marked Jacobian family, a contradiction.  Thus the generic
Artin invariant is two.  The locus $\sigma\leq1$ is closed
\cite[Corollary 10]{Ogu01}, giving the asserted nonempty open subset.
\end{proof}

\newcommand{\etalchar}[1]{$^{#1}$}

\Addresses
\end{document}